\documentclass[11pt]{amsart}

\usepackage{amsmath, amsthm, amsfonts}
\usepackage{epsfig} 
\newtheorem{prop}{Proposition}
\newtheorem{theorem}{Theorem}

\newtheorem{lemma}{Lemma}
\newtheorem{definition}{Definition}
\newtheorem{remark}{Remark}

\def\floor#1{\lfloor #1 \rfloor}
\def\ceil#1{\lceil #1 \rceil}
\newcommand{\Z}{\mathbb{Z}}
\newcommand{\N}{\mathbb{N}}

\newcommand{\un}{\mathbf{1}}
\renewcommand{\P}{\mathbb{P}}
\newcommand{\E}{\mathbb{E}}
\renewcommand{\S}{\mathcal{S}}
\newcommand{\F}{\mathcal{F}}
\newcommand{\G}{\mathcal{G}}

\author{Jean B\'{e}rard$^{1}$ and Alejandro Ram\'{\i}rez$^{1,2}$}

\thanks{$^1$Partially supported by ECOS-Conicyt grant CO5EO2}

\thanks{$^2$Partially supported by Fondo Nacional de Desarrollo Cient\'\i fico
y Tecnol\'ogico grant 1060738 and by Iniciativa Cient\'\i fica Milenio P-04-069-F}

\address[Jean B\'{e}rard]{\noindent Institut Camille Jordan, UMR CNRS 5208, 43, boulevard du 11 novembre
1918, Villeurbanne, F-69622, France; universit\'{e} de Lyon, Lyon, F-69003, France; 
universit\'{e} Lyon 1, Lyon, F-69003, France
\newline
e-mail:  \rm \texttt{jean.berard@univ-lyon1.fr}}

\address[Alejandro F. Ram\'\i rez]{Facultad de Matem\'aticas\\
Pontificia Universidad Cat\'olica de Chile\\
Vicu\~na Mackenna 4860, Macul\\
Santiago, Chile
\newline
e-mail:  \rm \texttt{aramirez@mat.puc.cl}}

\date{}

\title[CLT for the excited random walk in dimension $d \geq
2$]{Central limit theorem for the excited random walk in dimension $d \geq
2$}

\begin{document}

\begin{abstract}We prove that a law of large numbers and a central limit theorem hold for the excited random
  walk model in every dimension $d \geq
2$.\end{abstract}



\maketitle

\section{Introduction} An excited random walk with bias parameter $p\in(1/2,1]$
is a discrete time nearest neighbor random walk $(X_n)_{n \ge 0 }$ on the lattice $\Z^d$ obeying the 
following rule: when at time $n$ the walk is at a site  
 it has already visited before time $n$, it jumps uniformly at random to one of the $2d$ neighboring sites. 
 On the other hand, when the walk is at a site it has not visited before time $n$, it  
 jumps with probability $p/d$ to the right, probability $(1-p)/d$ to the left,
and probability $1/(2d)$ to the other
nearest neighbor sites. 

 The excited random walk was introduced in 2003 by 
Benjamini and Wilson~\cite{BenWil}, motivated by previous
works of~\cite{Dav, Dav2} and~\cite{PerWer} on self-interacting Brownian motions. 
Variations on this model have also been introduced.
The excited random walk on a tree was studied by Volkov~\cite{Vol}.
The so called {\it multi-excited random walk}, where the walk gets
pushed towards a specific direction upon its first $M_x$ visits to
a site $x$, with $M_x$ possibly being random, was introduced
by Zerner in~\cite{Zer1} (see also~\cite{Zer2} and~\cite{MouPimVal}).

In~\cite{BenWil}, Benjamini and Wilson proved
that for every value of $p\in (1/2,1]$ and $d \geq 2$, excited
random walks are transient. Furthermore, they proved
that for $d \geq 4$, 

\begin{equation}
\label{bw}
\liminf_{n\to\infty} n^{-1} X_n \cdot e_{1} >0 \quad a.s.,
\end{equation} where 
$(e_i:1\le i\le d)$ denote the canonical generators of the
group $\Z^d$.
Subsequently, Kozma extended (\ref{bw}) in~\cite{Koz1} and~\cite{Koz2}
to dimensions $d=3$ and $d=2$. Then, in~\cite{vdHHol}, relying on the lace expansion technique, van der Hofstad and Holmes  proved that a weak law of large 
numbers holds when $d>5$ and $p$ is close enough (depending on $d$) to $1/2$, and that a central limit theorem 
hold when $d>8$ and and $p$ is close enough (depending on $d$) to $1/2$.

In this paper, we prove that the biased coordinate of the
excited random walk satisfies a law of large numbers and
a central limit theorem for every $d \geq 2$ and $p \in (1/2,1]$.

\begin{theorem}\label{t:lgn-tcl} Let $p\in (1/2,1]$ and $d\ge 2$.

\begin{itemize}

\item[(i)] (Law of large numbers). There exists $v=v(p,d), \,      0<v<+\infty$ such that a.s.

$$
\lim_{n\to\infty}n^{-1} X_n\cdot e_1=v.
$$

\item[(ii)] (Central limit theorem). There exists $\sigma=\sigma(p,d), \, 
0<\sigma<+\infty$, such that
$$t \mapsto n^{-1/2}(X_{\floor{nt}}\cdot e_1 - v \floor{nt}),$$
converges in law as $n \to +\infty$ to a Brownian motion with 
variance $\sigma^2$.
\end{itemize}
\end{theorem}

Our proof is based on the well-known construction of regeneration times for the random walk, 
the key issue being to obtain good tail estimates for these regeneration times. Indeed, using estimates for the so-called
{\it tan} points of the simple random walk, introduced in~\cite{BenWil} and subsequently used in~\cite{Koz1, Koz2}, it is possible to prove that,
when $d \geq 2$, the number of distinct points visited by the excited random walk after $n$ steps is, with large probability, of order $n^{3/4}$ at least.
Since the excited random walk performs a biased random step at each time it visits a site it has not previously visited, the $e_{1}$-coordinate of the walk
 should typically be at least of order $n^{3/4}$ after $n$ steps. Since this number is $o(n)$, this estimate is not good enough to provide a direct proof 
 that the walk has linear speed. However, such an estimate is sufficient to prove that, while performing $n$ steps, 
 the walk must have many independent opportunities to perform a regeneration. A tail estimate on the regeneration times follows, and in turn, this  yields the law of large numbers and the central limit theorem, 
 allowing for a full use of the spatial homogeneity properties of the model.
 When $d \geq 3$, it is possible to replace, in our argument, estimates on the number of tan points by estimates on the 
 number of distinct points visited by the projection of the random walk on the $(e_{2},\ldots,e_{d})$ coordinates -- which is  essentially a simple random walk on $\Z^{d-1}$.
 Such an observation was used in~\cite{BenWil} to prove that~(\ref{bw}) holds when $d \geq 4$. 
 Plugging the estimates of~\cite{DonVar} in our argument, we can rederive the 
 law of large numbers and the central limit theorem when $d \geq 4$ without considering tan points. 
 Furthermore, a translation of the results in~\cite{Bol} and~\cite{Szn} about the volume of the 
 Wiener sausage to the random walk situation considered here, would allow us to  rederive our results when $d=3$, 
 and  to improve the tail estimates for any $d \geq 3$.

The regeneration time methods used to prove Theorem \ref{t:lgn-tcl} could also be used to describe the asymptotic behavior
of the configuration of the vertices as seen from the excited random walk. Let $\Xi:=\{0,1\}^{\Z^d \setminus \{0\}}$, equipped with the product topology and $\sigma-$algebra. 
For each time $n$ and site $x\ne X_n$, define $\beta(x,n):=1$ if the site $x$ was visited before time $n$ by the random walk, while $\beta(x,n):=0$
otherwise. Let $\zeta(x,n):=\beta(x-X_n,n)$ and define
$$\zeta(n):= (\zeta(x,n) ;  \,   x\in\Z^d \setminus \{0\} ) \in \Xi.$$
We call the process $(\zeta(n))_{n\in\N}$ the {\it environment seen from
the excited random walk}. 
 It is then possible to
show that if $\rho(n)$ is the law of $\zeta(n)$, there exists a probability measure $\rho$ defined on $\Xi$
such that

$$
\lim_{n\to\infty}\rho(n)=\rho,
$$
weakly.

In the following section of the paper we introduce the basic notation that 
will be used throughout. In Section \ref{rs}, we define the regeneration
times and formulate the key facts satisfied by them. In Section \ref{s:estimate} we obtain the tail estimates for the regeneration times via a good control
on the number of tan points. Finally, in Section \ref{sim}, we present the 
results of numerical simulations in dimension $d=2$ which suggest that, 
 as a function of the
bias parameter $p$, the speed $v(p,2)$ is an
increasing convex function of $p$, whereas the variance $\sigma(p,2)$
is a concave function which attains its maximum at some point strictly
between $1/2$ and $1$.

\section{Notations}
\label{not}

Let 
$\mathbf{b} := \{ e_{1},\ldots, e_{d}, -e_{1},\ldots, -e_{d} \}$.
Let $\mu$ be the distribution on $\mathbf{b}$ defined by
  $\mu(+e_{1})=p/d$, $\mu( -e_{1})=(1-p)/d$, $\mu(±\pm e_{j})=1/2d$ for $j \neq 1$.
  Let $\nu$ be the uniform distribution on $\mathbf{b}$.
Let $\S_{0}$ denote the sample space of the trajectories of the excited random walk starting at the origin:
$$    \S_{0}  := \left\{      (z_{i})_{i \geq 0} \in  (\Z^d)^{\N}; \,  z_{0}=0, \,  z_{i+1}-z_{i} \in \mathbf{b} \mbox{ for all $i \geq 0$}  \right\}.$$
 
For all $k \geq 0$, let $X_{k}$ denote the coordinate map defined on $\S_{0}$ by $X_{k}(  (z_{i})_{i \geq 0}   ) := z_{k}$.
We will sometimes use the notation $X$ to denote the sequence $(X_{k})_{k \geq 0}$.
We let $\F$ be the $\sigma-$algebra on $\S_{0}$ generated by the maps $(X_{k})_{k \geq 0}$.
For $k \in \N$, the sub-$\sigma-$algebra of $\F$ generated by $X_0,\ldots,X_k$ is denoted by $\F_k$.
And we let $\theta_{k}$ denote the transformation on
 $\S_{0}$ defined by
$(z_i)_{i \geq 0}   \mapsto  (z_{k+i} - z_k)_{i \geq 0} $.
For the sake of definiteness, we let $\theta_{+\infty}((z_i)_{i \geq 0}):=(z_i)_{i \geq 0}$.
For all $n \geq 0$, define the following two random variables on $(\S_{0},\F)$:
$$r_n := \max \{ X_i \cdot e_1 ; \, 0 \leq i \leq n \},$$ 
$$J_{n} = J_{n}(X) := \mbox{ number of indices $0 \leq k \leq n$ 
such that } X_{k} \notin \{ X_{i};    \,    0 \leq i \leq k-1 \}.$$
(Note that, with this definition, $J_{0}=1$.)

We now call $\P_0$ the law of  the excited random walk, which is formally defined as the unique probability measure on $(\S_{0},\F)$ 
 satisfying the following conditions: for every $k \geq 0$,
\begin{itemize}

\item on  $X_{k} \notin \{  X_{i}; \, 0 \leq i  \leq k-1 \}$, the conditional distribution of $X_{k+1}-X_{k}$ with respect to $\F_{k}$ is $\mu$;

\item  on  $X_{k} \in \{  X_{i}; \, 0 \leq i  \leq k-1 \}$, the conditional distribution of $X_{k+1}-X_{k}$ with respect to $\F_{k}$ is $\nu$.

\end{itemize}

\section{The renewal structure}
\label{rs}

We now define the regeneration times for the excited random walk
(see~\cite{SznZer} for the same definition 
in the context of random walks in random environment).
Define on $(\S_{0},\F)$ the following $(\mathcal F_k)_{k \geq 0}$-stopping times: $T(h) := \inf \{ k \geq 1; \, X_k \cdot e_1 > h \}$, and $D := \inf \{ k \geq 1; \, X_k \cdot e_1 = 0 \}$.
Then define recursively the sequences $(S_i)_{i \geq 0}$
and  $(D_i)_{i \geq 0}$ as follows:
$S_0:=T(0)$, $D_0 := S_{0} + D \circ \theta_{S_{0}}$, and
$S_{i+1} := T(r_{D_i})$, $D_{i+1} :=   S_{i+1} + D \circ \theta_{S_{i+1}}$ for $i\ge 0$,
with the convention that
$S_{i+1}=+\infty$ if $D_i = +\infty$, and, similarly, $D_{i+1}=+\infty$ if $S_{i+1} = +\infty$.
Then define 
$K := \inf \{ i \geq 0; \, D_i = +\infty \}$ and 
$\kappa := S_{K}$  (with the convention that $\kappa=+\infty$ when $K=+\infty$).

The key estimate for proving our results is stated in the following proposition. 

\begin{prop}\label{p:queue-regeneration}
As $n$ goes to infinity,
$$\P_0( \kappa \geq n ) \leq \exp \left( -n^{\textstyle{\frac{1}{19}  + o(1) } } \right).$$
\end{prop}

A consequence of the above proposition is that, under $\P_{0}$, $\kappa$ has finite 
moments of all orders, and also  $X_{\kappa}$, since the walk performs nearest-neighbor steps. 
We postpone the proof of 
Proposition~\ref{p:queue-regeneration} to Section~\ref{s:estimate}. 

\begin{lemma}\label{l:transience} There exists a $\delta>0$ such that $\P_0(  D  = +\infty)>\delta$.
\end{lemma}

\begin{proof}

This is a simple consequence of two facts. Firstly,  in~\cite{BenWil}
it is established that $\P_0$-a.s,
$ \lim_{k \to +\infty} X(k) \cdot e_1 = + \infty$.
On the other hand, a general lemma  (Lemma 9 of~\cite{Zer2}) shows that, given the first fact, an excited random walk satisfies 
$\P_0(  D  = +\infty)>0$.

\end{proof}

\begin{lemma}\label{l:sfini}
For all $h \geq 0$, $\P_{0}(T(h)<+\infty)=1$.
\end{lemma}
\begin{proof} This is
immediate from the fact that $\P_0$-a.s.,
$ \lim_{k \to +\infty} X(k) \cdot e_1 = + \infty$.
\end{proof}

Now define the sequence of regeneration times $(\kappa_n)_{n \geq 1}$ by
$\kappa_1:=\kappa$ and $\kappa_{n+1}:=\kappa_{n}+  \kappa \circ \theta_{\kappa_n} $,
with the convention that $\kappa_{n+1}=+\infty$ if $\kappa_{n}=+\infty$. 
For all $n \geq 0$, we denote by $\F_{\kappa_{n}}$ the completion with respect to $\P_{0}-$negligible sets 
of the $\sigma-$algebra generated by the events of the form
$\{ \kappa_{n} = t  \} \cap A$, for all $t \in \N$, and $A \in \F_{t}$.

The following two propositions are analogous respectively to Theorem 1.4 and Corollary 1.5  of~\cite{SznZer}. Given Lemma~\ref{l:transience} and Lemma~\ref{l:sfini}, 
the proofs are completely similar to those presented in~\cite{SznZer},
noting that the process $(\beta(n),X_n)_{n\in\N}$ is strongly Markov, so we omit them, and refer the reader to~\cite{SznZer}. 

\begin{prop}\label{p:renouvellement}

For every $n \geq 1$, $ \P_{0}(\kappa_{n}<+\infty) = 1$. Moreover, for every $A \in \F$, 
the following equality holds $\P_{0}-$a.s.
\begin{equation}\label{e:renouveau}\P_{0} \left(  X \circ \theta_{\kappa_{n}}  \in  A | \F_{\kappa_{n}}  \right) = \P_{0} \left(  X  \in  A |  D = +\infty  \right).\end{equation}

\end{prop}

\begin{prop}\label{p:renouvellement3}

With respect to $\P_0$, the random variables $\kappa_1$, $\kappa_2 - \kappa_1$, $\kappa_3 - \kappa_2,\ldots$
are independent, and, for all $k \geq 1$, the distribution of $\kappa_{k+1} - \kappa_k$ with respect to $\P_0$ is that 
of $\kappa$ with respect to $\P_0$ conditional upon $D=+\infty$. Similarly, 
 the random variables $X_{\kappa_1}$, $X_{\kappa_2} - X_{\kappa_1}$, $X_{\kappa_3} - X_{\kappa_2},\ldots$
are independent, and, for all $k \geq 1$, the distribution of $X_{\kappa_{k+1}} - X_{\kappa_k}$ with respect to $\P_0$ is that 
of $X_{\kappa}$ with respect to $\P_0$ conditional upon $D=+\infty$.

\end{prop}

For future reference, we state the following result.

\begin{lemma}\label{l:temps-arret}

On $S_{k}<+\infty$, the conditional distribution of the sequence
$(X_{i} - X_{S_{k}})_{S_k  \leq i < D_k }$ with respect to $\F_{S_{k}}$ is the same as the distribution 
of $(X_{i})_{0 \leq i < D}$ with respect to $\P_0$. 
\end{lemma}

\begin{proof}
Observe that  between times $S_{k}$ and $D_{k}$, the walk never visits any site that it has visited before time $S_{k}$. Therefore, applying the
strong Markov property to the process $(\beta(n),X_n)_{n\in\N}$ and spatial translation
invariance, we conclude the proof.
\end{proof}

A consequence of Proposition~\ref{p:queue-regeneration} is that $\E_0( \kappa | D=+\infty) < +\infty$ and   
$\E_0( |X_{\kappa}| | D=+\infty) < +\infty$. Since $\P_{0}(\kappa \geq 1)=1$ and $ \P_{0}(X_{\kappa} \cdot e_{1} \geq 1)=1$, 
$\E_0( \kappa | D=+\infty)>0$ and $\E_0( X_{\kappa} \cdot e_{1} | D=+\infty) >0$. Letting  $v(p,d):= \frac{\E_0( X_{\kappa} \cdot e_{1}| D=+\infty)}{  \E_0( \kappa | D=+\infty)  }$,
we see that  $0 < v(p,d)<+\infty$.

The following law of large numbers can then be proved, using Proposition~\ref{p:renouvellement3}, exactly as 
Proposition 2.1 in~\cite{SznZer}, to which we refer for the proof. 

\begin{theorem}\label{t:lgn}
Under $\P_{0}$, the following limit holds almost surely:
$$\lim_{n \to +\infty} n^{-1} X_{n} \cdot e_{1} =v(p,d).$$
\end{theorem}

Another consequence of Proposition~\ref{p:queue-regeneration} is that $\E_0( \kappa^{2} | D=+\infty) < +\infty$ and   
$\E_0( |X_{\kappa}|^{2} | D=+\infty) < +\infty$. Letting 
$\sigma^{2}(p,d):= \frac{\E_0(  \left[X_{\kappa} \cdot e_{1}   -  v(p,d) \kappa \right]^{2} | D=+\infty)}{  \E_0( \kappa | D=+\infty)  }$,
we see that $\sigma(p,d)<+\infty$. That $\sigma(p,d)>0$ is explained in Remark~\ref{r:sigmapositif} below.

The following functional central limit theorem can then be proved, using Proposition~\ref{p:renouvellement3}, exactly as 
Theorem 4.1 in~\cite{Szn2}, to which we refer for the proof. 

\begin{theorem}\label{t:tcl}
Under $\P_{0}$, the following convergence in distribution holds: as $n$ goes to infinity,
$$t \mapsto n^{-1/2}(X_{\floor{nt}}\cdot e_1 - v \floor{nt}),$$
converges to a Brownian motion with 
variance $\sigma^2(p,d)$.
\end{theorem}
\begin{remark}\label{r:sigmapositif}
The fact that $\sigma(p,d)>0$ is easy to check. Indeed, we will
prove that the probability of the event $X_{\kappa}\cdot e_{1} \ne v \kappa$
is positive.
There is a positive probability that the first step of the walk is $+e_{1}$, and that $X_{n} \cdot e_{1} > 1$ for all $n$ afterwards. 
In this situation, $\kappa=1$ and $X_{\kappa} \cdot e_{1} = 1$. Now, there is a positive probability that the walk first performs 
the following sequence of steps: $+e_{2}, -e_{2}, +e_{1}$, and that then $X_{n} \cdot e_{1} > 1$ for all $n$ afterwards. In this situation, 
$\kappa=3$ and $X_{\kappa}\cdot e_{1} = 1$. 
\end{remark}

\section{Estimate on the tail of $\kappa$}\label{s:estimate}

\subsection{Coupling with a simple random walk and tan points}

We use the coupling of the excited random walk with a simple random walk that was 
introduced in~\cite{BenWil}, and subsequently
used in~\cite{Koz1, Koz2}.

To define this coupling, let 
$(\alpha_{i})_{i \geq 1}$ be a sequence of i.i.d. random variables 
with uniform distribution on the set $\{ 1, \ldots, d \}$.
Let also $(U_{i})_{i \geq 1}$ be an i.i.d. family of random variables with uniform
distribution on $[0,1]$, independent from $(\alpha_{i})_{i \geq 1}$.  
Call $(\Omega,\G,P)$ the probability space on which these variables are defined.
Define the sequences of random variables 
$Y=(Y_{i})_{i \geq 0}$ and $Z=(Z_{i})_{i \geq 0}$ taking values in $\Z^{d}$, as follows. First, 
$Y_{0}:=0$ and $Z_{0}:=0$. Then consider $n \geq 0$, and assume that $Y_{0},\ldots,Y_{n}$ and $Z_{0},\ldots,Z_{n}$ have 
already been defined.
Let $Z_{n+1}:= Z_{n} +      ( \un ( U_{n+1} \leq 1/2)   -  \un ( U_{n+1} > 1/2)   ) e_{\alpha_{n+1}} $. 
Then, if $Y_{n} \in \{ Y_{i};    \,    0 \leq i \leq n-1 \}$ or $\alpha_{n+1} \neq 1$, 
let  $Y_{n+1}:= Y_{n} +      ( \un ( U_{n+1} \leq 1/2)   -  \un ( U_{n+1} > 1/2)   ) e_{\alpha_{n+1}} $. 
Otherwise, let $Y_{n+1} := Y_{n} +      ( \un ( U_{n+1} \leq p)   -  \un ( U_{n+1} > p)   ) e_{1}$.

The following properties are then immediate:

\begin{itemize}
\item $(Z_{i})_{i \geq 0}$ is a simple random walk on $\Z^{d}$;
\item $(Y_{i})_{i \geq 0}$ is an excited random walk on $\Z^{d}$ with bias parameter $p$;
\item for all $2 \leq j \leq d$ and $i \geq 0$,  $Y_{i} \cdot e_{j} = Z_{i} \cdot e_{j}$;
\item the sequence $( Y_{i} \cdot e_{1} - Z_{i} \cdot e_{1} )_{i \geq 0}$ is non-decreasing.
\end{itemize}

\begin{definition}
If $(z_{i})_{i \geq 0} \in \S_{0}$,  we call an integer $n \geq 0$ an $(e_{1},e_{2})$--tan point index for the sequence 
$(z_{i})_{i \geq 0}$ 
if $z_{n} \cdot e_{1} > z_{k} \cdot e_{1}$ for
all $0 \leq  k \leq n-1$ such that $z_{n} \cdot e_{2} = z_{k} \cdot e_{2}$. 
\end{definition}

The key observation made in~\cite{BenWil} is the following. 

\begin{lemma}\label{l:tan-point}
If $n$ is an $(e_{1},e_{2})$--tan point index for  $(Z_{i})_{i \geq 0}$, 
then $Y_{n} \notin  \{ Y_{i};    \,    0 \leq i \leq n-1 \}$. 
\end{lemma}

\begin{proof}
If $n$ is an $(e_{1},e_{2})$--tan point index and if there exists an $\ell \in \{  0, \ldots, n-1 \}$ such that
$Y_{n}=Y_{\ell}$, then observe that,  using the fact that 
$Z_{\ell} \cdot e_{2} = Y_{\ell} \cdot e_{2}$ and $Z_{n} \cdot e_{2} = Y_{n} \cdot e_{2}$,
we have that $Z_\ell\cdot e_2=Z_n\cdot e_2$. Hence, by the definition of a tan point
we must have that $Z_{\ell} \cdot e_{1} < Z_{n} \cdot e_{1}$,  whence
$Y_{n}\cdot e_{1} - Z_{n} \cdot e_{1} <  Y_{\ell}\cdot e_{1} - Z_{\ell} \cdot e_{1}$.
But this contradicts the fact that the coupling has the property that
 $Y_{n}\cdot e_{1} - Z_{n} \cdot e_{1} \geq  Y_{\ell}\cdot e_{1} - Z_{\ell} \cdot e_{1}$.

\end{proof}

Let $H := \{   i \geq 1 ;  \,   \alpha_{i} \in \{  1, 2  \}       \}$, and define the sequence of indices $(I_{i})_{i \geq 0}$
by $I_{0}:=0$, $I_{0} < I_{1} < I_{2} < \cdots$, and $ \{  I_{1},I_{2},\ldots  \} = H $. 
Then the sequence of random variables $(W_{i})_{i \geq 0}$ defined by 
$W_{i}:= (Z_{I_{i}}  \cdot e_{1} , Z_{I_{i}} \cdot e_{2} )$ 
forms a simple random walk on $\Z^{2}$. 

If $i$ and $n$ are such that $I_{i}=n$, it is immediate to check that $n$ is an $(e_{1},e_{2})$--tan point index
for $(Z_{k})_{k \geq 0}$ if and only if $i$ is an $(e_{1},e_{2})$--tan point index for the random walk $(W_{k})_{k \geq 0}$.

For all $n \geq 1$, let $N_{n}$ denote the number of $(e_{1},e_{2})$--tan point indices of $(W_{k})_{k \geq 0}$
that are $\leq n$.
The arguments used to prove the following lemma are quite similar to the 
ones used in the proofs of Theorem 4 in~\cite{BenWil} and Lemma 1 in~\cite{Koz2}, which are themselves partly based
on estimates in~\cite{BouSch}.

\begin{lemma}\label{l:queue-tan}
For all $0 < a < 3/4$,  as $n$ goes to infinity, 
$$P(N_{n} \leq n^{a}) \leq \exp \left( - n^{  \textstyle{\frac{1}{3}-\frac{4a}{9}} + o(1)   } \right).$$
 \end{lemma}

\begin{proof}

For all $k \in \Z \setminus \{  0  \}$, $m \geq 1$, consider the three sets
\begin{eqnarray*} \Gamma(m)_{k} &:=&  \Z \times \{   2k \floor{ m^{1/2}}   \},\\
\Delta(m)_{k} &:=& \Z \times (   (2k-1)    \floor{m^{1/2}}  ,    (2k+1)    \floor{m^{1/2}}      ) ,\\
\Theta(m)_{k} &:=& \{  v \in \Delta(m)_{k}; \,    |  v   \cdot e_{2} | \geq    2k \floor{m^{1/2}}      \}.\end{eqnarray*}

Let $\chi(m)_{k}$ be the (a.s. finite since the simple random walk on $\Z^{2}$ is recurrent) first time when $(W_{i})_{i \geq 0}$ 
hits $\Gamma(m)_{k}$. Let $\phi(m)_{k}$ be the (again a.s. finite for the same reason as $\chi(m)_{k}$) 
first time after $\chi(m)_{k}$ when $(W_{i})_{i \geq 0}$ leaves $\Delta(m)_{k}$. 
Let $M_k(m)$ denote the number of time indices $n$ that are $(e_{1},e_{2})$--tan point indices,  
and satisfy $ \chi(m)_{k}   \leq  n  \leq \phi(m)_{k}-1$ and 
$W_{n} \in \Theta(m)_{k}$. 
 
Two key observations in~\cite{BenWil}  (see Lemma 2 in~\cite{BenWil} and the discussion before its statement) are that: 
\begin{itemize}
\item the sequence $(M_k(m))_{k \in \Z\setminus \{ 0 \}}$ is i.i.d.;
\item there exist $c_{1},c_{2}>0$ such that $P(M_1(m) \geq c_{1} m^{3/4} ) \geq c_{2}$.
\end{itemize}

Now, consider an $\epsilon>0$ such that $b:=1/3-4a/9-\epsilon>0$. 
let $m_{n}:= \ceil{   \left(  n^{a}/c_{1} \right)^{4/3}}+1$, 
and let $h_{n} := 2  ( \ceil{n^{b}}  + 1  )  \floor{m_{n}^{1/2}} $. 
Note that, as $n \to +\infty$,  $(h_n)^{2} \sim  (4 c_{1}^{-4/3})  n^{\frac{2}{3}+\frac{4}{9}a-2\epsilon}$.
Let $R_{n,+}$ and $R_{n,-}$ denote the following events
$$R_{n,+} := \{  \mbox{for all } k \in \{1,\ldots,  + \ceil{n^{b}}   \}, \, M_{k}(m_{n}) \leq    c_{1}  m_{n}^{3/4}       \},$$
and
$$R_{n,-} := \{  \mbox{for all } k \in \{-\ceil{n^{b}},\ldots,-1\}, \, M_{k}(m_{n}) \leq    c_{1}  m_{n}^{3/4}       \}.$$
From the above observations, $P(R_{n,+} \cup R_{n,-}) \leq 2 (1-c_{2})^{\ceil{n^{b}}}$.

Let $q_{n}:=\floor{n (h_{n})^{-2}}$, and let $V_{n}$ be the event 
$$V_{n} :=   \{   \mbox{for all } i \in \{0,\ldots,n\} , \,   -h_{n} \leq  W_{i} \cdot e_{2} \leq +h_{n}            \}  .$$
By Lemma~\ref{l:sortie} below, there exists a constant $c_{3}>0$ such that, for all large enough $n$, all $ -h_{n}      \leq  y   \leq +h_{n}$,  and $x \in \Z$, 
the probability that a simple random walk on $\Z^{2}$ started at $(x,y)$ at time zero leaves   $ \Z \times \{-h_{n},\ldots,  +h_{n}\}$ before time 
$h_{n}^{2}$, is larger than $c_{3}$. 
A consequence is that, for all $q \geq 0$, the probability that the same walk fails to leave  $ \Z \times \{-h_{n},\ldots,  +h_{n}\}$ before time 
$q h_{n}^{2}$ is less than $(1-c_{3})^{q}$. Therefore 
$P(V_{n}) \leq (1-c_{3})^{q_{n}}$. 

Observe now that, on $V_{n}^{c}$,  
$$n \geq \max(\phi(m_{n})_{k} ; \, 1 \leq k \leq \ceil{n^{b}} ) \mbox{ or }  n\geq \max(\phi(m_{n})_{k} ; \, -  \ceil{n^{b}} \leq k \leq -1 ).$$
Hence, on $V_{n}^{c}$, 

$$N_{n} \geq \sum_{k=1}^{ \ceil{n^{b}}} M_k(m_{n}) \mbox{ or } N_n\geq \sum_{k=-1}^{ - \ceil{n^{b}} } M_k(m_{n}).$$ 
We deduce that, on $R_{n,+}^{c} \cap R_{n,-}^{c}  \cap    V_{n}^{c}$, $N_{n} \geq c_{1} m_{n}^{3/4} > n^{a}$.

As a consequence, $\P_{0}(N_{n}  \leq n^{a}) \leq \P_{0}(R_{n,+} \cup R_{n,-}) + \P_{0}(V_{n})$, 
so that 
$\P_{0}(N_{n}  \leq n^{a}) \leq 2 (1-c_{2})^{\ceil{n^{b}}} + (1-c_{3})^{q_{n}}$

Noting that, as $n$ goes to infinity,  $q_{n} \sim n h_{n}^{-2} \sim    (4 c_{1}^{-4/3})^{-1}  n^{1/3-4a/9+2 \epsilon}$,
the conclusion follows.

\end{proof}

\begin{lemma}\label{l:sortie}
There exists a constant $c_{3}>0$ such that, for all large enough $h$, all $ -h      \leq  y   \leq +h$,  and $x \in \Z$, 
the probability that a simple random walk on $\Z^{2}$ started at $(x,y)$ at time zero leaves   $ \Z \times \{-h,\ldots,  +h\}$ before time 
$h^{2}$, is larger than $c_{3}$. 

\end{lemma}

\begin{proof}

Consider the probability that the $e_{2}$ coordinate is larger than $h$ at time $h^{2}$. By standard coupling, this probability is minimal
when $y=-h$, so the central limit theorem applied to the walk starting with $y=-h$ yields the existence of $c_{3}$.

\end{proof}

%

\begin{lemma}\label{l:tan-dimension-d}
For all $0 < a < 3/4$,  as $n$ goes to infinity, 
$$\mathbb{P}_{0}(J_{n} \leq n^{a}) \leq \exp \left( - n^{\textstyle{\frac{1}{3}-\frac{4a}{9}} + o(1) } \right).$$
\end{lemma}

\begin{proof}

Observe that, by definition, $I_{k}$ is the sum of $k$ i.i.d. random variables whose distribution is geometric with parameter $2/d$.
By a standard large deviations bound, there is a constant $c_{6}$ such that, for all large enough $n$,  
$P(I_{\floor{n d^{-1}}} \geq n) \leq \exp(-c_{6} n)$.
Then observe that, if $I_{\floor{n d^{-1}}} \leq n$,  we have
$J_{n}(Y) \geq N_{\floor{n d^{-1}}}$ according to Lemma~\ref{l:tan-point} above.
(Remember that, by definition, $J_{n}(Y)$ is the number of indices $0 \leq k \leq n$ such that  
$Y_{k} \notin \{ Y_{i};    \,    0 \leq i \leq k-1 \}$. )
Now, according to Lemma~\ref{l:queue-tan} above, we have that,
for all $0 < a < 3/4$,  as $n$ goes to infinity, 
$$P(N_{\floor{n d^{-1}}} \leq \floor{nd^{-1}}^{a}) \leq \exp \left( -   \floor{nd^{-1}}^{ \textstyle{\frac{1}{3}-\frac{4a}{9}} + o(1)  } \right),$$
from which it is easy to deduce that for all  $0 < a < 3/4$,  as $n$ goes to infinity, 
$P(N_{\floor{n d^{-1}}} \leq n^{a}) \leq \exp \left( - n^{  \textstyle{\frac{1}{3}-\frac{4a}{9}} + o(1)  } \right)$.
 Now we deduce from the union bound that
 $P(  J_{n} (Y)  \leq n^{a} ) \leq P(I_{\floor{n d^{-1}}} \geq n) +  P(   N_{\floor{n d^{-1}}} \leq n^{a})$.
 The conclusion follows.

\end{proof}

\subsection{Estimates on the displacement of the walk}

\begin{lemma}\label{l:borne-deplacement}

For all $1/2 < a < 3/4$,  
as $n$ goes to infinity,
$$ \P_0(  X_n \cdot e_{1} \leq  n^{a})  \leq   \exp \left( - n^{ \psi(a) + o(1)  } \right),$$
where   $\psi(a):=  \min \left( \textstyle{\frac{1}{3}-\frac{4a}{9}},      2a - 1 \right)$.

\end{lemma}

\begin{proof}
Let $\gamma:= \frac{2p-1}{2d}$.
Let $(\varepsilon_{i})_{i \geq 1}$ be an i.i.d. family of random variables with common distribution $\mu$ on $\mathbf{b}$, and
let  $(\eta_{i})_{i \geq 1}$ be an i.i.d. family of random variables with common distribution $\nu$ on  $\mathbf{b}$
 independent from $(\varepsilon_{i})_{i \geq 1}$. Let us call $(\Omega_{2},\G_{2},Q)$ the probability space on which these variables are defined.

Define the sequence of random variables  $(\xi_{i})_{i \geq 0}$ taking values in $\Z^{d}$, as follows.
First, set $\xi_{0}:=0$.
 Consider then $n \geq 0$, assume that $\xi_{0},\ldots,\xi_{n}$
have 
already been defined, and consider the number 
$J_{n}(\xi)$ of indices $0 \leq k \leq n$ such that  
$\xi_{k} \notin \{ \xi_{i};    \,    0 \leq i \leq k-1 \}$. 
If $\xi_{n} \notin \{ \xi_{i};    \,    0 \leq i \leq n-1 \}$, set 
$\xi_{n+1} := \xi_{n} + \varepsilon_{J_{n}(\xi)}$. 
Otherwise, let $\xi_{n+1} := \xi_{n} + \eta_{n  -  J_{n}(\xi)+1}$.
It is easy to check that the sequence $(\xi_{n})_{n \geq 0}$ is an excited random walk on $\Z^{d}$ with bias parameter $p$.

Now, according to Lemma~\ref{l:tan-dimension-d}, 
for all $1/2<a<3/4$, $Q(J_{n} \leq  n^{a}) \leq \exp \left( - n^{  \textstyle{\frac{1}{3}-\frac{4a}{9}} + o(1)   } \right)$.
That for all $1/2<a<3/4$, $Q(J_{n-1} \leq 2 \gamma^{-1} n^{a}) \leq \exp \left( - n^{  \textstyle{\frac{1}{3}-\frac{4a}{9}} + o(1)  } \right)$
is an easy consequence. Now observe that, by definition, for $n \geq 1$,
$\xi_{n} = \sum_{i=1}^{J_{n-1}}  \varepsilon_{i}    +  \sum_{i=1}^{n-J_{n-1}}  \eta_{i}$. 
Now, there exists a constant $c_{4}$ such that, for all large enough $n$, and every $  2 \gamma^{-1} n^{a}   \leq k \leq n$, 
$$Q \left(     \sum_{i=1}^{k}     \varepsilon_{i}\cdot e_{1}  \leq  (3/2) n^{a}   \right) \leq Q\left(     \sum_{i=1}^{k}     \varepsilon_{i}\cdot e_{1}  \leq \textstyle{\frac{3}{4} \gamma k}     \right) \leq \exp \left( -c_{4} n^{a}          \right) ,$$
by a standard large deviations bound for the sum $ \sum_{i=1}^{k}     \varepsilon_{i}\cdot e_{1}$, whose terms are i.i.d. bounded random variables with expectation $\gamma>0$.
By the union bound, we see that 
$$Q\left(     \sum_{i=1}^{J_{n-1}}     \varepsilon_{i}\cdot e_{1}  \leq        (3/2) n^{a} \right) \leq  n \exp \left( -c_{4} n^{a}          \right)        +\exp \left( - n^{  \textstyle{\frac{1}{3}-\frac{4a}{9}} + o(1)  } \right).$$

Now, there exists a constant $c_{5}$ such that, for all large enough $n$, and 
for every $1 \leq k \leq n$, 
$Q \left(     \sum_{i=1}^{k}     \eta_{i}\cdot e_{1}  \leq  -(1/2) n^{a}   \right) \leq \exp \left( -c_{5} n^{2a-1} \right) ,$
by a standard Gaussian upper bound for the simple symmetric random walk on $\Z$.

By the union bound again, we see that 
$Q \left(     \sum_{i=1}^{n-J_{n-1}}   \eta_{i}\cdot e_{1}  \leq  -(1/2) n^{a}   \right) \leq n \exp \left( -c_{5} n^{2a-1} \right)$.
The conclusion follows.

\end{proof}

\begin{lemma}\label{l:borne-temps-de-retour}

As $n$ goes to infinity, 
$$ \P_0( n \leq D < +\infty )  
\leq \exp \left( - n^{ 1/11+ o(1)} \right).$$ 

\end{lemma}

\begin{proof}
Consider $1/2<a<3/4$, and write $\P_0( n \leq D < +\infty ) = \sum_{k=n}^{+\infty}  \P_0(D = k) \leq \sum_{k=n}^{+\infty}  \P_0(X_{k} \cdot e_{1} = 0)  \leq  \sum_{k=n}^{+\infty}  \P_0(X_{k} \cdot e_{1} \leq k^{a}).$
Now, according to Lemma~\ref{l:borne-deplacement},
$ \P_0(  X_k \cdot e_{1} \leq  k^{a})  \leq   \exp \left( - k^{ \psi(a) + o(1)  } \right)$.
It is then easily checked that 
$\sum_{k=n}^{+\infty}   \exp \left( - k^{ \psi(a) + o(1)  } \right) \leq \exp \left( - n^{ \psi(a) + o(1)  } \right)$.
As a consequence, $\P_0( n \leq D < +\infty ) \leq \exp \left( - n^{ \psi(a) + o(1)  } \right)$.
Choosing $a$ so as to minimize $\psi(a)$, the result follows.

\end{proof}

\subsection{Proof of Proposition~\ref{p:queue-regeneration}}

Let $a_{1}, a_{2}, a_{3}$ be positive real numbers such that $a_{1}<3/4$ and $a_{2}+a_{3}<a_{1}$. 
For every $n>0$,  let $u_n:=\floor{n^{a_{1}}}$, $v_n :=  \floor{n^{a_{2}}}$,  $w_n :=  \floor{n^{a_{3}}}$.
In the sequel, we assume that $n$ is large enough so that  $v_{n} (w_{n}+1) +2 \leq u_{n}$. 
Let $$A_n := \{  X_n  \cdot e_{1} \leq  u_{n} \}; \,    B_{n} := \bigcap_{k=0}^{v_{n}} \{ D_k < +\infty  \}; \,  C_{n}:=  \bigcup_{k=0}^{v_{n}} \{ w_{n} \leq D_k - S_k < +\infty \}.$$
(With the convention that, in the definition of $C_{n}$, $D_{k}-S_{k}=+\infty$ whenever $D_{k}=+\infty$.)
We shall prove that 
$\{ \kappa \geq n   \} \subset A_{n} \cup B_{n} \cup C_{n}$, then apply the union bound to 
$\P_{0}(    A_{n} \cup B_{n} \cup C_{n}        )$, and then separately bound the three probabilities
$\P_{0}(A_{n})$, $\P_{0}(B_{n})$, $\P_{0}(C_{n})$.

 Assume that $A_n^c \cap B_n^c \cap C_n^c $ occurs. Our goal is to prove that this assumption implies that $\kappa < n$. 

Call $M$ the smallest index $k$ between $0$ and $v_{n}$ such that 
$D_k=+\infty$, whose existence is ensured by $B_n^c$. By definition, $\kappa=S_{M}$, so we have to prove that $S_{M} < n$. 
For notational convenience, let $D_{-1}=0$. 
By definition of $M$, we know that $D_{M-1}<+\infty$. 
Now write $r_{D_{M-1}}=\sum_{k=0}^{M-1}   (r_{D_{k}} - r_{S_{k}})    + (r_{S_{k}} - r_{D_{k-1}})$, with the convention that
$\sum_{k=0}^{-1}=0$.
Since the walk performs nearest-neighbor steps, we see that for all $0 \leq k \leq M-1$, $r_{D_{k}} - r_{S_{k}} \leq D_{k} - S_{k}$.
On the other hand, by definition, for all $0 \leq k \leq M-1$, $r_{S_{k}} - r_{D_{k-1}} =1 $.
Now,  for all $0 \leq k \leq M-1$, 
$D_k-S_k \leq w_{n}$, due to the fact that $C_{n}^{c}$ holds and that $D_{k}<+\infty$. 
As a consequence, we obtain that $r_{D_{M-1}}   \leq  M  (w_{n}+1) \leq  v_{n} (w_{n}+1)$.
Remember now that $v_{n} (w_{n}+1) +2 \leq u_{n}$, so we have proved that, 
$r_{D_{M-1}}+2 \leq u_{n}$. Now observe that, on $A_{n^{c}}$, $X_{n} \cdot e_{1} > u_{n}$. As a consequence, the smallest $i$ such that $X_{i} \cdot e_{1} =  r_{D_{M-1}}+1$
must be $< n$. But  $S_{M}$ is indeed the smallest $i$ such that $X_{i} \cdot e_{1} =r_{D_{M-1}}+1$, so we have proved that $S_{M}<n$ on $A_{n}^{c} \cap B_{n}^{c} \cap C_{n}^{c}$.

The union bound then yields the fact that, for large enough $n$, 
 $\P_{0}(\kappa \geq n) \leq \P_{0}(A_{n})+ \P_{0}(B_{n})+ \P_{0}(C_{n})$.

Now, from Lemma~\ref{l:borne-deplacement}, we see that
$ \P_0(  A_n )  \leq  \exp( -n^{\psi(a_{1}) + o(1)} ).$ 
By repeatedly applying Lemma~\ref{l:sfini} and the strong Markov property at the stopping times $S_k$ for $k=0,\ldots, v_{n}$ to the
process $(\beta(n),X_n)_{n\in\N}$, we see that 
$\P_0(B_n) \leq \P_0(D<+\infty)^{v_{n}}$. Hence, from Lemma~\ref{l:transience}, we know that $\P_0(B_n)\leq (1-\delta)^{v_n}$.

From  the union bound and Lemma~\ref{l:temps-arret}, we see that 
$\P_{0}(C_n) \leq (v_{n}+1) \P_0( w_{n} \leq D < +\infty )$, so, by Lemma~\ref{l:borne-temps-de-retour},
$\P_{0}(C_n) \leq  (v_{n}+1) \exp( - n^{a_{3}/11    + o(1)})$.

Using Lemma~\ref{l:borne-deplacement},  we finally obtain the following estimate:
$$\P_{0}(\kappa \geq n) \leq  (1-\delta)^{\floor{n^{a_{2}}}} +  (  \floor{n^{a_{2}}}   +1)\exp\left( - n^{a_{3}/11    + o(1)}\right)  + \exp \left( -n^{\psi(a_{1}) + o(1)} \right).$$

Now, for all $\epsilon$ small enough, choose $a_{1}=12/19$, $a_{2}=1/19$, $a_{3}=11/19-\epsilon$. This ends the proof of Proposition~\ref{p:queue-regeneration}.

\section{Simulation results}
\label{sim}

We have performed simulations of the model in dimension $d=2$, using a C code and the Gnu Scientific Library
facilities for random number generation. 

The following graph is a plot of an estimate of $v(p,2)$ as a function of $p$.
Each point is the average over 1000 independent simulations of $(X_{10000} \cdot e_{1})/10000$.

\includegraphics*[width=7cm, angle=+90]{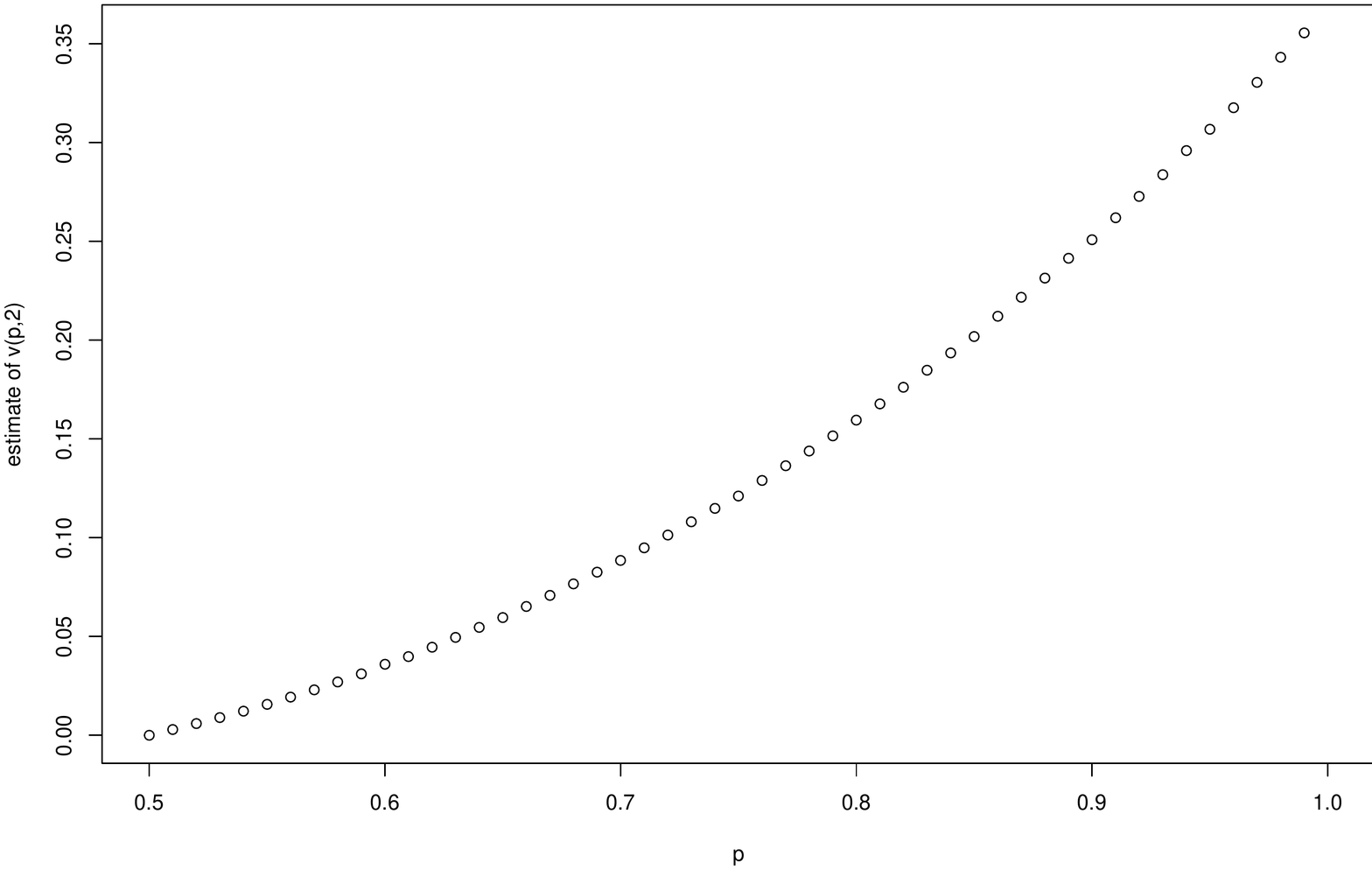}

The following graph is a plot of an estimate of $\sigma(p,2)$ as a function of $p$.
Each point is the standard deviation over 1000000 independent simulations 
of $(X_{10000} \cdot e_{1})/(10000)^{1/2}$ (obtaining a reasonably smooth curve required many more simulations for
$\sigma$ than for $v$).

\includegraphics*[width=7cm, angle=+90]{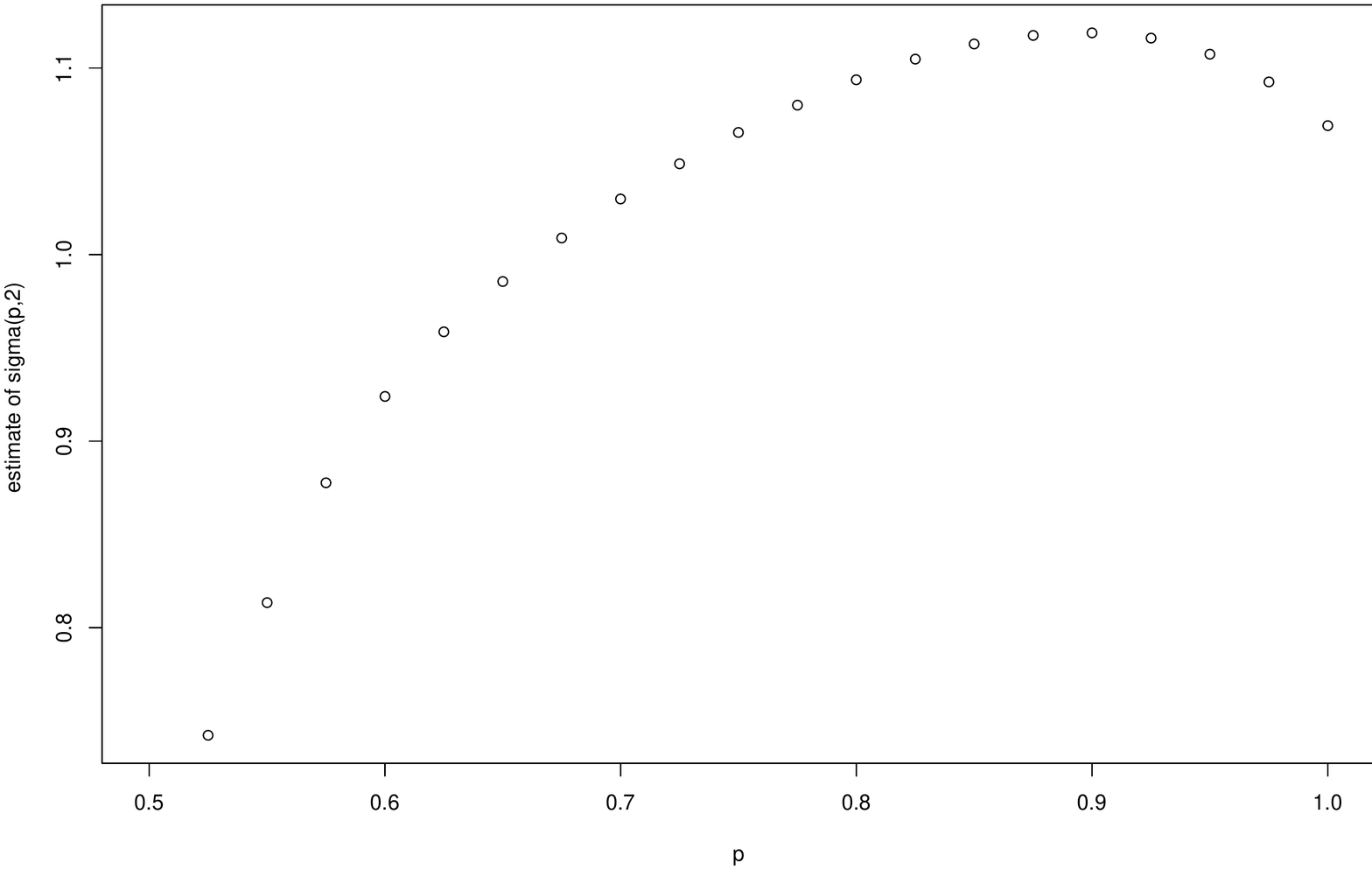}

\bibliographystyle{plain}
\bibliography{clt-excited-rw-article}

\end{document}